\documentclass[12pt,a4paper]{article}
\usepackage{xypic,amssymb,amsthm,a4wide}
\newtheorem{thm}{Theorem}[section]
\newtheorem{corol}[thm]{Corollary}
\newtheorem{lemma}[thm]{Lemma}
\newtheorem{prop}[thm]{Proposition}
\newtheorem{defin}[thm]{Definition}
\newtheorem{rem}[thm]{Remark}
\newtheorem{ex}[thm]{Example}
\newenvironment{example}{\begin{ex}\rm}
       {\hfill$\vartriangle$\end{ex}}
\newenvironment{remark}{\begin{rem}\rm}
       {\hfill$\vartriangle$\end{rem}}
\let\cal=\mathcal\let\frak=\mathfrak\let\bb=\mathbb
\newcommand{\Cc}{{\cal C}}
\newcommand{\E}{{\cal E}}
\renewcommand{\L}{{\cal L}}
\newcommand{\F}{{\cal F}}

\newcommand{\Pcc}{{\cal P}}

\newcommand{\R}{{\bb R}}
\newcommand{\Z}{{\bb Z}}
\newcommand{\C}{{\bb C}}

\newcommand{\X}{{\hat X}}
\newcommand{\T}{{\hat T}}
\newcommand{\fp}{{X\times_B \X}}

\newcommand{\Hom}{\mbox{\rm Hom}}
\newcommand{\Alt}{\mbox{\rm Alt}}

\newcommand{\Pcic}{\mbox{\rm Pic}}
\newcommand{\obj}{\mbox{\rm Ob}}
\newcommand{\Mor}{\mbox{\rm Mor}}
\newcommand{\loc}{\mathbf{Loc}(\T)}
\newcommand{\sky}{{\mathbf{Sky}(T)}}
\newcommand{\bysame}{$\raise.2em\hbox to 3em{\hrulefill}$\thinspace, }
\begin{document}
\thispagestyle{empty}
\rightline{To be published (in two parts) in \emph{J. Geom. Phys.}} \vfill
\begin{center}
{\bfseries\Large
A Fourier transform for sheaves  on Lagrangian  \\[10pt]
families of real tori }
\par\addvspace{20pt}
{\sc U. Bruzzo,$^{\S\P}$  G. Marelli$^{\S}$ \ {\rm and} \ F. Pioli$^{\S}$}
\par\medskip
$^{\S}$\thinspace Scuola Internazionale Superiore di Studi Avanzati (SISSA),\\
Via Beirut 4, 34013 Trieste, Italy \\
$^{\P}$\thinspace Istituto Nazionale di Fisica Nucleare, \\ Sezione di Trieste
\end{center}
\vfill
\begin{center} 24 May 2001\\Revised 15 June and 18 July 2001 \end{center}
\vfill
\begin{quote} \small {\sc Abstract.}
We systematically develop a transform of the Fourier-Mukai type
for sheaves on symplectic manifolds $X$ of any dimension fibred in Lagrangian
tori. One obtains a bijective correspondence between unitary local systems
supported on Lagrangian submanifolds of $X$ and holomorphic  vector
bundles with compatible unitary connections supported on complex submanifolds
of the relative Jacobian of $X$ (suitable conditions being verified on both
sides).   \end{quote}
\vfill
\leftline{\hbox to8cm{\hrulefill}}\par
{\footnotesize
\noindent This work has been partly supported by  the Friuli-Venezia Giulia
Region Research
Project ``Noncommutative geometry: algebraic, analytical and probabilistic
aspects and applications
to mathematical physics" and by  a grant awarded to the third author by the
University of Genova
under the programme ``Finanziamenti a singoli e giovani ricercatori." The
first and third author
are  members of VBAC (Vector Bundles on Algebraic Curves), which is
partially supported by {\sc Eager}
(EC FP5 Contract no. HPRN-CT-2000-00099) and  {\sc Edge} (EC FP5 Contract
no. HPRN-CT-2000-00101).\par
\par
\noindent\emph{E-Mail addresses:} {\tt bruzzo@sissa.it, marelli@sissa.it,
pioli@sissa.it}
}
\eject

\section{Introduction} The idea that, in accordance with the
Strominger-Yau-Zaslow conjecture
\cite{SYZ}, a kind of Fourier-Mukai transform
should describe transformation properties of D-branes under string-theoretic
mirror symmetry dates back to 1996 \cite{DM}. The original Fourier-Mukai
transform,
mapping coherent sheaves
on an abelian variety $X$ to coherent sheaves on the dual variety
$\X$, was introduced in \cite{M2}. A relative Fourier-Mukai transform
for elliptic varieties was developed in \cite{BBHP1,B,HP,BBHP2}
and was shown to  describe a correct D-brane transformation
pattern in the case of K3 surfaces \cite{BBHP1,AD}.
An analogous result
was shown to hold for elliptic Calabi-Yau threefolds in \cite{ACHY}.

In the case of Calabi-Yau threefolds  which are fibred in (special
Lagrangian) real 3-tori, a
similar description should be provided by a ``real'' relative Fourier
transform. The presence of
singular fibres raises here a big problem because it is not clear how to
handle them. As a first step, one may consider the simplified case when
there are no singular fibres.

If $X$ is a symplectic family of smooth Lagrangian  tori,
the dual family $\X$ has a natural complex structure.
Then the relative   Fourier transform yields a
correspondence between local systems supported on Lagrangian
submanifolds of $X$ and holomorphic  vector
bundles supported on complex submanifolds of $\X$,
where both sets of data satisfy suitable conditions (cf.~later in this
introduction). Some results along
these lines were already contained in \cite{AP} but we strengthen and
extend them
considerably (cf. also \cite{LYZ}). We also carefully spell out
the conditions on the  submanifold $S$ of $X$ which ensure
that the support of the transformed sheaf is a complex submanifold of $\X$.

The correspondence we get closely resembles Fukaya's homological
mirror symmetry \cite{Fu}. Comparison with that approach suggests
that in order to extend the results presented in this paper
to more general Lagrangian submanifolds (e.g.,  when the Lagrangian
submanifold $S$ is ramified
over the base of the fibration $X$), or to the situation when $X$ has
singular fibres,
it is necessary to allow for some kind of quantum corrections.
Future extensions of the theory should also
allow for the inclusion of a B-field and should investigate
the possibility of describing the correspondence between the
Floer homology of $X$ and the algebraic cohomology of $\X$ in terms
of the Fourier-Mukai transform studied in this paper.
One could also study the relation between the transform presented
in this paper and the constructions of Laumon \cite{L} and Rothstein
\cite{Roth}; the local systems we transform are $\cal D$-modules,
the same objects considered by these authors. This might also
relate to possible applications to generalizations of the
Krichever correspondence along the lines of \cite{Nak}.  

We describe now the contents of this paper.
As a preliminary step, in section \ref{absolute}
we consider the absolute case. In its simplest version
one constructs two functors which establish the
equivalence between the category $\sky$ of skyscraper sheaves
of finite-dimensional vector spaces on a real torus $T$,
and the category $\loc$ of unitary local systems (vector bundles with a flat
unitary connection) on the dual torus $\T$. We also study the
transformation of a local system supported on an affine subtorus of $T$.

In section \ref{relative} we consider the relative case. In \ref{framework}
we set up the general
framework. Then, in considering local systems supported on a Lagrangian
submanifold of a symplectic torus fibration $X\to B$, we first analyze the
two extreme cases, i.e., when the submanifold is a fibre of $X$ or a
Lagrangian section   (section \ref{fibres+sections}). In the first case one
gets the usual tautological property of the Fourier-Mukai transforms, while in
the second one obtains  a bijective  correspondence between local systems
supported on Lagrangian sections of $X$ and holomorphic  bundles with
compatible unitary connections, flat along the fibre directions of $\X$ (and
satisfying some further conditions).

The intermediate, non-transversal cases (i.e., when one considers a Lagrangian
submanifold $S\subset X$ whose projection onto $B$
has a dimension strictly between 0 and $\dim B$) are more
involved, and are analyzed in section
\ref{skew}. To get a well-behaved transform one needs
to assume that $S$ intersects the fibres
$X_b$ of $X$ (here $b\in B$) along subtori $S_b$ of $X_b$, and  that
the vertical tangent spaces to $S$   undergo parallel displacement
under the natural Gauss-Manin connection defined in $TX$. Under this
condition the transform of a local system on $S$ is a holomorphic vector bundle
supported on a complex submanifold of $\X$, and again,  provided that suitable
conditions are satisfied, there is a bijective correspondence between the two
sets of data.

These result hold true
whatever is the dimension of $X$, and do not require $X$ to be Calabi-Yau
(and not
even complex). One should note that,
when $X$ is a Calabi-Yau manifold,
the additional condition on the support $S$ we have previously described is
in general quite unrelated to the
condition of
$S$ being special (in addition to being Lagrangian), and coincides with the
latter only
when $X$ is complex 3-dimensional, and
the projection of the Lagrangian submanifold onto the base is (real)
1-dimensional
(this corresponds to a transformed sheaf which is a line bundle supported
on a curve in $\X$).
It is not clear to the authors whether this situation has any implication
or motivation in string theory.

In section \ref{Conclusions} we draw some conclusions, in particular
we comment upon the relation of this construction to Fukaya's homological
mirror
symmetry.

Two Appendices contain two proofs that, sketchy as they are, are too
lengthy to be
included in the main text.

\par\smallskip
{\bf Acknowledgements.} We are   thankful to K.~Fukaya and D.~Hern\'an\-dez
Ruip\'erez for useful discussions.

\section{The absolute case\label{absolute}}
We start this section by  offering a description of (smooth) $U(1)$  bundles
on real tori in terms of their factor of automorphy  which fully parallels the
one available (in the holomorphic case) on complex tori (cf.~\cite{LB}). This
description of line bundles will be extensively used in the remainder of the
paper. In section \ref{complexes} we describe
two complexes which are naturally associated with the
Poincar\'e sheaf. In section \ref{equivalence}  we introduce and study
the functors which establish the equivalence between the category of
skyscrapers and that of local systems; the computation of the second functor
will require the study of the cohomology of a complex associated
with the Poincar\'e bundle.
In section \ref{subtori} we consider local systems supported on subtori.
\subsection{Line bundles on real tori\label{bundles}} Let $\Lambda$ be
a $g$-dimensional lattice in a $g$-dimensional
real vector space $V$, and let $T=V/\Lambda$ be
the corresponding torus. Let $\Pcic(T)$ denote
the group of isomorphism classes of $U(1)$ bundles on $T$.
The group $\Pcic(T)$ is isomorphic to a group
$P(\Lambda)$ we may associate with the lattice $\Lambda$
in the following way. As a set, $P(\Lambda)$ is the
set of pairs $(A,\chi)$, where $A\in\Alt^2(\Lambda,\Z)$
is an alternating two-form on $\Lambda$, and
$\chi$ is a \emph{semicharacter} for $A$, namely,
a map $\chi\colon\Lambda\to U(1)$ such that
$$\chi(\lambda+\mu)=\chi(\lambda)\,\chi(\mu)\,e^{i\pi A(\lambda,\mu)}$$
for all $\lambda,\mu\in\Lambda$. The group structure
is the one given by
$$(A_1,\chi_1)\cdot (A_2,\chi_2)=(A_1+A_2,\chi_1\,\chi_2).$$
The isomorphism $\Pcic(T)\simeq P(\Lambda)$ is the Appell-Humbert theorem
for real tori. Via the isomorphism $\Alt^2(\Lambda,\Z)\simeq H^2(T,\Z)$,
the form $A$ is to be identified with the first Chern class.
In this way we have an exact sequence
$$ 0 \to \Hom_\Z(\Lambda,U(1)) \to \Pcic(T) \stackrel{c_1}{\relbar\joinrel\to}
H^2(T,\Z) \to 0 $$
and the kernel $\Hom_\Z(\Lambda,U(1))$, whose elements
are isomorphism classes of flat $U(1)$ bundles, is isomorphic
to the dual torus $\T=V^\vee/\Lambda^\vee$ (here $V^\vee=\Hom_\R(V,\R)$,
$\Lambda^\vee=\Hom_\Z(\Lambda,\Z)$). To every point $y\in \T$
there corresponds a flat line bundle $\L_y$ whose associated
pair is
$$A_y=0,\qquad \chi_y(\lambda)=e^{2i\pi\,y(\lambda)}.$$

The description of the bundle $\L$ by means of the pair $(A,\chi)$
allows one to give an explicit characterization of the global
sections of $\L$. To this end one introduces
the \emph{factor of automorphy} of the  pair $(A,\chi)$, defined as the map
$$a_\L\colon V\times\Lambda\to U(1),\qquad a_\L(x,\lambda)=
\chi(\lambda)\,e^{i\pi A(x,\lambda)}$$
(here $A$ has been extended to $V\times V$ in the natural way).

\begin{prop} Let $\L$ be a line bundle on $T$, corresponding to the
pair $(A,\chi)\in P(\Lambda)$. The global sections of $\L$ are in
a one-to-one correspondence with the smooth functions
$s\colon V\to\C$ satisfying the automorphy condition
$$s(x+\lambda)=a_\L(x,\lambda)\,s(x)$$
for all $x\in V$, $\lambda\in\Lambda$.
\end{prop}
\begin{proof} The proof is a (simplified) replica of the one
holding in the case of complex tori \cite{LB}.\end{proof}

The action of an automorphism of $\L$ changes the factor of automorphy; an
automorphism of $\L$ is induced by a map $\phi\colon V\to U(1)$, and the new
factor of automorphy is
$$a'_\L(x,\lambda)=\phi(x+\lambda)\,a_\L(x,\lambda)\,\phi(x)^{-1}.$$

Now we use these tools to describe the
\emph{Poincar\'e bundle}
$\Pcc$ on the product $T\times\T$. The line bundle
$\Pcc$ is associated with the pair $(A,\chi)\in P(\Lambda\times\Lambda^\vee)$,
where
$$A((\lambda_1,\mu_1),(\lambda_2,\mu_2)=\mu_1(\lambda_2)-\mu_2(\lambda_1),\qquad
\chi(\lambda,\mu)= e^{i\pi\mu(\lambda)}.$$
The corresponding factor of automorphy is
$$a_{\Pcc}(x,y,\lambda,\mu)=e^{i\pi[y(\lambda)-\mu(x)-\mu(\lambda)]}.$$
It is convenient to apply the automorphism induced by the map
$$\phi\colon V\times V^\vee\to U(1),\qquad \phi(x,y)=e^{i\pi\,y(x)}$$
thus obtaining a new factor of automorphy
\begin{equation}\label{aut}
a'_{\Pcc}(x,y,\lambda,\mu)=e^{2i\pi\,y(\lambda)}.
\end{equation}
This description of the Poincar\'e bundle shows explicitly
that $\Pcc_{\vert T\times\{y\}}\simeq \L_y$.

The connection form of the connection $\nabla_\Pcc$ on the Poincar\'e bundle
is written in the gauge where the factor of
automorphy of $\Pcc$ has the form (\ref{aut}) as
\begin{equation}\label{con}\mathbb{A} = 2i\,\pi\sum_{j=1}^gy_j\,dx^j
\end{equation}
where $(x^1,\dots,x^g)$ are flat coordinates on $T$ and
$(y_1,\dots,y_g)$ are dual flat coordinates on $\hat T$.

If we act on $a_{\Pcc}$ with the automorphism $\phi(x,y)=e^{-i\pi\,y(x)}$
we obtain the factor of automorphy
$a''_{\Pcc}(x,y,\lambda,\mu)=e^{-2i\pi\,\mu(x)}$
which shows that, after the
identification $\hat{\kern-1.7pt\mbox{$\hat T$}}\simeq T$, the dual bundle
$\Pcc^\vee$ is a Poincar\'e bundle for $\T\times T$.

\subsection{Complexes associated with the Poincar\'e sheaf\label{complexes}}
We denote by $p $,
$\hat p $ the projections onto the two factors of $T\times\T$.
To simplify notation we shall set\footnote{A word about notation:
if $f\colon X\to Y$ is a differentiable map between two differentiable
manifolds, and $\cal F$ any sheaf on $Y$, we shall denote by $f^{-1}\cal F$
the sheaf-theoretic
inverse  image of $\cal F$; if $\cal F$ is a sheaf of $\Cc_Y^\infty$-modules,
we shall denote by $f^\ast\cal F$ its inverse image as a sheaf of modules,
i.e.,
$$f^\ast\cal F=f^{-1}\cal F\otimes_{f^{-1}\Cc_Y^\infty}\Cc_X^\infty\,.$$}
$$\Omega^{m,n}=p ^\ast\Omega^m_T\otimes_{\Cc^\infty_{T\times\T}}
\hat p^\ast\Omega^n_\T\,.$$
The  connection $\nabla_\Pcc$  has a K\"unneth splitting
into two operators
$$\nabla_1\colon\Pcc\to\Pcc\otimes\Omega^{1,0},\qquad
\nabla_2\colon\Pcc\to\Pcc\otimes\Omega^{0,1}$$
both squaring to zero (but their anticommutator is the
curvature of $\nabla_\Pcc$).  In the gauge of equation (\ref{aut}), the action
of $\nabla_1$, $\nabla_2$ on sections is locally written in the form
\begin{equation}\label{nablas}\nabla_1s=\sum_{j=1}^g\left(\frac{\partial
s}{\partial x^j}+2i\pi\,y_j\,s\right)\,dx^j,\qquad
\nabla_2f=\sum_{j=1}^g\frac{\partial f}{\partial
y_j}\,dy_j,\end{equation}
\noindent where $g$ is the dimension of $T$.

Let $\E$ be a $\Cc^\infty_T$-module with a flat connection $\nabla$.
By pulling the pair $(\E,\nabla)$ back to $T\times\T$ and coupling it
with the pair $(\Pcc,\nabla_1)$ we obtain
a complex
$$\begin{array}{lcl}
0 &\to& \ker\nabla_1^{\E} \to p^\ast\E\otimes\Pcc
\stackrel{\nabla_1^{\E}}{\relbar\joinrel\to} p^\ast\E\otimes\Pcc \otimes
\Omega^{1,0} \\
&\stackrel{\nabla_1^{\E}}{\relbar\joinrel\to} &
p^\ast\E\otimes\Pcc \otimes
\Omega^{2,0} \to \dots \end{array}$$
Since locally the operator $\nabla_1^{\E}$ coincides with the
exterior differential, this sheaf complex is exact, and is a fine
resolution of the sheaf $\ker\nabla_1^{\E}$. Thus we obtain an
isomorphism
$$H^i(T\times U, \ker\nabla_1^{\E}) \simeq
H^i((p^\ast\E\otimes\Pcc \otimes
\Omega^{\bullet,0})(T\times U),\nabla_1^{\E}),\qquad i\ge 0$$
between the  cohomology of the sheaf $\ker\nabla_1^{\E}$
and the cohomology of the complex $(p^\ast\E\otimes\Pcc \otimes
\Omega^{\bullet,0})(T\times U)$ (where $U$ is an open set in $\T$)
acted upon by the differential $\nabla_1^{\E}$. The sheaf $R^i\hat
p_\ast\ker\nabla_1^{\E}$
associated with the presheaf $U\rightsquigarrow H^i(T\times U,
\ker\nabla_1^{\E})$
is the $i$th higher direct image of the sheaf $\ker\nabla_1^{\E}$ to $\T$.

The same results hold for the operator $\nabla_2$.
The cohomology of the complex $(\Gamma(\Pcc^\vee \otimes
\Omega^{0,\bullet}),\nabla_2)$ will be computed in
section \ref{equivalence}.

\subsection{The equivalence \label{equivalence}}
For every real torus $T$ we consider

(i) the category $\sky$ of skyscrapers on $T$ of total finite length
(i.e., $\dim H^0(T,M)<\infty$ for all $M\in\obj(\sky))$;

(ii) the category $\mathbf{Loc}(T)$ of unitary local systems on $T$.
Objects of this category are pairs $(E,\nabla)$, where
$E$ is a smooth complex vector bundle on $T$, and $\nabla$
is a flat unitary connection. Morphisms in this category
are vector bundle morphisms compatible with the connections.
The objects in $\mathbf{Loc}(T)$ can also be regarded
as locally free $\C_T$-modules of finite rank equipped with
a hermitian metric defined up to homothety; under this disguise,
a unitary local system will be typically denoted by a gothic letter. We shall
freely switch from one picture to the other.

We now define the functor $\F\colon\sky\to\loc$ \cite{AP}.
Let $M$ be a skyscraper of finite length on $T$
supported on a 0-dimensional subset $S$, with
the $C^\infty_T$-structure given by evaluation
of functions. Denoting by $\Pcc_S$ the restriction
of $\Pcc$ to $S\times T$, and by $p_S$, $\hat p_{ S}$
the projections of $S\times \T$ onto its factors,
we define $\hat\E$ as the sheaf $\hat p_{S,\ast}(p_S^\ast M
\otimes \Pcc_S)$. This
is locally free of finite rank, so that it is the sheaf of sections of a
vector bundle
$\hat E$, with $\mbox{rk}\,\hat E=\mbox{length}(M)$. Moreover,
the operator $\nabla_2$ naturally extends to $p^\ast_{S}  M\otimes\Pcc_S$,
and induces an operator
$\hat\nabla:\hat \E\to \hat \E\otimes\Omega^1_\T$ which is a unitary flat
connection.
Standard checks show that this procedure does define  a functor.

\begin{example} Let $\C(x)$ denote the one-dimensional
skyscraper at $x\in T$. One has $\F(\C(0))\simeq\C_\T$.
Indeed, in this case we have $p^\ast\mathcal M\otimes\Pcc\simeq
\Cc^\infty_{\{0\}\times\T}$, and, in view of equations
(\ref{con}), (\ref{nablas}), the operator $\nabla_2$ reduces on this sheaf
to the exterior differential along the $\T$ direction. As
a consequence, $(\hat E,\hat\nabla)=(\Cc^\infty_\T,d)$, and
$\F(\C(0))=\ker d \simeq\C_\T$.
\label{ex1}\end{example}

For every $x\in T$ let $t_x$ be the associated translation,
$t_x(x')=x+x'$. Moreover, identify $\hat{\kern-1.7pt\mbox{$\hat T$}}$ with
$T$. The following result is easily proved.
\begin{prop} For every $x\in T$ and
$M\in\obj(\sky)$ there is an isomorphism $\F(t_x^{-1} M)\simeq
L_{-x}\otimes\F(M)$.
\label{trasl}\end{prop}
As a consequence, in view of Example  \ref{ex1}, we have
\begin{corol} For every $x\in T$ one has
$\F(\C(x))\simeq L_{-x}$.
\label{trasf}\end{corol}
This defines the action of the functor $\F$ on the whole
category $\sky$.

Now we construct the inverse functor.
Let $(\hat E,\hat \nabla)$ be  an object in $\loc$,  and let $\hat \E$ be the
sheaf of sections of $\hat E$. As we did in the previous section,
but reversing the roles of $T$ and $\T$, we consider
on the sheaf $\hat p^\ast\hat \E\otimes\Pcc^\vee$
an operator $\hat \nabla_2^{\hat \E}$  obtained by coupling
(the pullback of) $\hat \nabla$ with the operator $\hat \nabla_2$
(the $\hat T$-component of the connection of $\Pcc^\vee$).
We shall eventually prove the following.
\begin{prop} 1.
$R^jp_\ast\ker\hat \nabla_2^{\hat \E}=0$ for $j=0,\dots,g-1$;

2. The sheaf $R^gp_\ast\ker\hat \nabla_2^{\hat \E}$ is  a skyscraper of
finite length.
\end{prop}
The functor $\hat \F$ is defined as $\hat \F((\hat E,\hat \nabla))=
R^gp_\ast\ker\hat \nabla_2^{\hat \E}$.

As a first step we compute the action of $\hat\F$ on the
trivial line bundle, i.e, we take $\hat E=\Cc^\infty_\T$ and
$\hat \nabla=d$.   Thus we want to compute the sheaves
$R^jp_\ast\ker\hat \nabla_2$. To this end we shall study the presheaves
$$U\rightsquigarrow H^j(U\times\T,\ker\hat \nabla_2)\simeq
H^j\left((\Pcc^\vee\otimes\Omega^{0,\bullet})(U\times\T),\hat \nabla_2\right)$$
whose associated sheaves are exactly the sheaves we are interested in.

As a first result we have
\begin{prop} $H^0(U\times\T,\ker\hat \nabla_2)=0$ for all open sets
$U\subset T$, so that $p_\ast\ker\hat \nabla_2=0.$
\label{vanish}\end{prop}
\begin{proof}
An element of $H^0(U\times\T,\ker\hat \nabla_2)$
restricted to $\{x\}\times\T$, with $x\in U$, yields a global
section of $L_x$, which is zero unless $x=0$. By a density
argument we get the result.
\end{proof}

\noindent To compute the higher-order direct images
we first consider the case $g=1$.
\begin{prop} \label{prop1} If $g=1$, then  $R^1p_\ast\ker\hat\nabla_2\simeq
\C(0)$.
\label{R1}\end{prop}
\begin{proof} We compute the cohomology groups $H^1(U\times\T,\ker\hat\nabla_2)
\simeq H^1((\Pcc^\vee\otimes\Omega^{0,\bullet})(U\times\T),\hat\nabla_2)$.
We represent $T$ as
$\R=/\Z\lambda$, with
$\lambda\in\R^\ast$, and $\T=\R/\Z\mu$, with $\mu=1/\lambda$. Let $W$ be
the inverse image of $U$ in $\R$.

We work now in a gauge where the factor of automorphy of
$\Pcc^\vee$ is $e^{2i\pi\,\mu(x)}$, and the operator
$\hat\nabla_2$ is the $\T$-part of the exterior differential.
An element in $((\Pcc^\vee\otimes\Omega^{0,1})(U\times\T),\ker\hat\nabla_2)$
may be written as $\tau= t(x,y)\,dy$, where $t$ is a function on
$W\times V^\vee$  satisfying the
automorphy condition
$$t(x,y+\mu)=t(x,y)\,e^{2i\pi\,\mu(x)}.$$
If $\tau$ is a coboundary, $\tau=\hat\nabla_2s$, one has
$$s(x,y)=\int_0^y\,t(x,u)\,du+c(x).$$
The function $s$ must satisfy the automorphy condition,
which amounts to the following condition on $c$:
\begin{equation}\label{cond}
c(x)(1-e^{2i\pi\,\mu(x)})=-\int_0^\mu\,t(x,u)\,du.
\end{equation}
If $0\notin U$ this condition may be solved for $c$, so that
$H^1(U\times\T,\ker\hat\nabla_2)=0$. Thus
$R^1p_\ast\ker\hat\nabla_2$ is a skyscraper supported at $0\in T$.

If $0\in U$, the condition (\ref{cond}) may be solved if and only
if $$\int_0^\mu\,t(0,u)\,du=0$$
so that $H^1(U\times\T,\ker\hat\nabla_2)\simeq \C$. This proves
the claim.
\end{proof}

We move to the higher-dimensional case by means of a K\"unneth-type
argument.

\begin{prop}\label{prop2} If $\dim T=g$ we have

1. $R^jp_\ast\ker\hat\nabla_2=0$ for $j=0,\dots,g-1$;

2. $R^gp_\ast\ker\hat\nabla_2\simeq \C(0)$.
\end{prop}
\begin{proof}
A choice of flat coordinates $(x^1,\dots,x^g)$ on $T$
fixes an isomorphism $T\simeq S^1\times\dots S^1$.
The Poincar\'e sheaf $\Pcc$ on $T\times\T$ is the product
of the Poincar\'e sheaves $\Pcc_i$ on the $i$ factors of $T\times\T$, as
one can check for instance by describing the
Poincar\'e bundles by their factors of automorphy.
Let
$U\subset T$ be of the form $U=U_1\times\dots\times U_g$
where each $U_i$ lies in a factor of $\T$. If $g=2$, a word-by-word
translation of the K\"unneth theorem for de Rham cohomology
(cf.~e.g.~\cite{BT}) gives a decomposition
$$H^j(U\times\T,\ker\hat\nabla_2)=\bigoplus_{m+n=j}
H^m(U_1\times S^1,\ker\hat\nabla_2^1)\otimes H^n(U_2\times
S^1,\ker\hat\nabla_2^2)$$
whence we have, by Proposition \ref{vanish},
$$H^j(U\times\T,\ker\hat\nabla_2)=0\quad\mbox{for}\quad j=0,1,\quad
H^2(U\times\T,\ker\hat\nabla_2)\simeq \C.$$
Induction on $g$ then yields, for every $g$,
$$H^j(U\times\T,\ker\hat\nabla_2)=0\quad\mbox{for}\quad j=0,\dots,g-1,\quad
H^g(U\times\T,\ker\hat\nabla_2)\simeq \C.$$
This proves both claims.
\end{proof}
So we have also obtained
$$H^j(T\times\T,\ker\hat\nabla_2)=
\left\{\begin{array}{ccl} 0 & \mbox{for} & j=0,\dots,g-1, \\
\C & \mbox{for} & j=g.\end{array}\right.$$

Let $\mathfrak L_x$ be the local system corresponding to the
line bundle $\L_x$ with its flat connection. In analogy with Proposition
\ref{trasl}, we have

\begin{prop} $\hat \F(\mathfrak L_{-x}\otimes_{\C_\T}\mathfrak S)\simeq
t_x^{-1}\hat\F(\mathfrak S)$ for every $x\in T$ and every local system
$\mathfrak S$ on
$\T$.
\end{prop}

\begin{corol} $\hat\F(\mathfrak L_{-x})\simeq \C(x)$ for every $x\in T$.
\label{complete}\end{corol}

\begin{remark} Since any flat vector bundle on a torus is a direct sum of
flat line bundles (i.e., every local system on $\T$
is a direct sum of local systems of the type $\mathfrak L_x$),
Corollary \ref{complete} completely describes the action of the
functor $\hat\F$.
\label{remark}\end{remark}

Corollaries \ref{trasf} and \ref{complete} and Remark \ref{remark}
eventually prove
\begin{thm} The functors $\F$, $\hat\F$ are inverse to each other, and
establish
an equivalence between the categories $\sky$ and $\loc$.\label{thm}
\end{thm}
Again, any question related to the behaviour of morphisms under the functor
$\hat\F$ is simply a matter of routine checks.

\subsection{Subtori \label{subtori}}
We consider now the transformation of $U(1)$ local
systems supported on affine subtori of the $g$-dimensional torus $T=V/\Lambda$.

\begin{defin} A subtorus of $T$ is a subset $S\subset T$ of the form
$S=W/W\cap\Lambda$, where $W$ is $k$-dimensional linear subspace
of $V$ such that the lattice $W\cap\Lambda$ has rank $k$. An affine
subtorus is a subset of the form $S+x$ for an element $x\in T$.
\end{defin}

Let $\frak L\equiv (\L,\nabla)$ be a $U(1)$ local system supported on
a $k$-dimensional affine subtorus $S$ of $T$.
By coupling the pullback of $\nabla$ with the connection of
$\Pcc_S$ and projecting on the $T$-differentials
one obtains a complex
$$\begin{array}{lcl}
0 &\to& \ker\nabla_1^{\cal L} \to p^\ast_S\cal L \otimes\Pcc_S
\stackrel{\nabla_1^{\cal L}}{\relbar\joinrel\to} p^\ast_S \cal L
\otimes\Pcc_S \otimes
\Omega^{1,0}  \\
&\stackrel{\nabla_1^\cal L}{\relbar\joinrel\to} &
p^\ast_S\cal L \otimes\Pcc_S \otimes
\Omega^{2,0}  \to \dots \end{array}$$
\begin{prop}\label{absoluteprop}
\begin{enumerate}\item  $R^j\hat p_{S,\ast} \ker\nabla_1^{\cal L}=0$ for all
$j\ne k$; \item  $R^k\hat p_{S,\ast}  \ker\nabla_1^{\cal L}$
is supported on a $(g-k)$-dimensional affine subtorus $\hat S$ of $\hat T$,
which is normal to $S$;
\item if $\cal L$ is trivial, then $\hat S$ goes through the origin
of $\hat T$, otherwise it is an affine subtorus translated
by the element of $\hat T$ corresponding to $\cal L^\ast$.
\item The sheaf  $R^k\hat p_{S,\ast} \ker\nabla_1^{\cal L}$ on $\hat S$ is
a $U(1)$ bundle, and has a
compatible flat  connection which makes it into a $U(1)$ local system
$\hat\frak L$.
\end{enumerate}\end{prop}
\begin{proof} The proof of this Proposition in given in   Appendix A.
\end{proof}
Let us describe the content of this Proposition in local coordinates; while
this
is just simple linear algebra, the explicit equations we are going to write
will
help to understand the more complicated relative situation.
Let $(y^1,\dots,y^g)$ be flat coordinates in $T$, $(w_1,\dots,w_g)$
the corresponding dual flat coordinates in the dual torus $\hat T$,
and write the equation for the affine subtorus
$S$ in the form
$$\sum_{j=1}^g a_j^i\,y^j+\chi^i=0,\qquad i=1,\dots,g-k\,.$$
The equations  $\sum_{j=1}^g a_j^i\,y^j=0$ describe a corresponding
``linear subtorus''
$S_0$; the equations of the dual torus $S_0^\ast$ may be written implicitly
as
$$\sum_{j,\ell=1}^g a_j^i\,g^{j\ell}\,w_\ell=0\,,\qquad i=1,\dots,g-k\,,$$
where the constant functions $g^{j\ell}$ are the components of the natural flat
metric on  $\hat T$, or explicitly as
\begin{equation} w_\ell=\sum_{m=1}^k\tilde
a^m_\ell\,\xi_m\,,\qquad\ell=1,\dots,g
\label{explicit}\end{equation}
for a suitable $k\times(g-k)$ matrix $\tilde a$.
The specification of the local system $\frak L$ corresponds to a choice of
the parameters
$(\xi_1,\dots,\xi_k)$ in equation (\ref{explicit}).
The support $\hat S$ of the transformed local system is given by equations
$$\sum_{j=1}^g \gamma_{m}^j\,w_j+\xi_m=0\,,\qquad m=1,\dots,k\,,$$
where $\gamma_{\ell}^j$ is a matrix satisfying $\sum_{j=1}^g \gamma_{\ell}^j\,a_j^i=0$.
The local system $\hat\frak L$ is given by the point in $\hat S_0^\ast$  whose coordinates are the
numbers $\chi^i$.

The pair $(\hat S,\hat\frak L)$ is the \emph{Fourier-Mukai
transform} of the pair $(S,\frak L)$.
Of course we may perform the same transformation from $\hat T$ to $T$ (in
addition
to the obvious replacements, one twists by $\Pcc^\vee$ instead of $\Pcc$),
and we have:
\begin{prop} The Fourier-Mukai transform of $(\hat S,\hat\frak L)$
is naturally isomorphic to the pair $(S,\frak L)$.
\end{prop}

Let  $\mathbf{Loc}_k(T)$ be the category of $U(1)$ local systems
supported on
affine subtori of $T$ of dimension $k$. Objects of this category are triples
$(S,\L,\nabla)$ (where $S$ is an affine subtorus in $T$, $\L$  a line
bundle on $S$, and $\nabla$ a flat unitary connection on $\L$) modulo
isomorphisms, i.e., modulo vector bundle isomorphisms which commute with the
actions of the connections (the two line bundles having the same support).
The space of  morphisms between two objects
$(S_1, \L_1, \nabla_1)$ and $(S_2, \L_2, \nabla_2)$ of $\mathbf{Loc}_k(T)$  is
defined by taking into account that the  intersection $S= S_1\cap
S_2$ is a (possibly empty) finite collection of (possibly zero-dimensional)
affine tori $R_i$, and one sets $$\Mor ((S_1, \L_1, \nabla_1), (S_2, \L_2,
\nabla_2)) =  \bigoplus_i \Mor_{\nabla} ((R_i, \L_1, \nabla_1), (R_i, \L_2,
\nabla_2)), $$ where $\Mor_{\nabla} ((R_i, \L_1, \nabla_1), (R_i, \L_2,
\nabla_2))$ is  the set of morphisms between ${\L_1}_{\vert R_i}$ and
${\L_2}_{\vert R_i}$ compatible with the connections $\nabla_1 $ and $\nabla_2
$. It is easy to check that the Fourier-Mukai transform yields
an equivalence of categories $$ \mathbf{Loc}_k(T)\simeq
\mathbf{Loc}_{g-k}(\hat T)\,.$$

\section{Relative theory\label{relative}}
\subsection{The geometric setting\label{framework}}
Let $(X,\omega)$ be a connected  symplectic manifold admitting
a map $f\colon X\to B$ whose fibres
are $g$-dimensional smooth Lagrangian tori. We assume that $f$
admits a Lagrangian section $\sigma\colon B\to X$;
according to \cite{Dui}, this makes $X$ isomorphic, as a symplectic
manifold fibred in Lagrangian submanifolds,
to a quotient bundle $T^\ast B/\Lambda$, where
$\Lambda$ is a Lagrangian covering of $B$. The symplectic form $\omega$
provides an isomorphism $\mbox{Vert}\,{TX}\simeq f^\ast T^\ast B$.
We also have an identification
$T B\simeq R^1f_\ast\R\otimes\C^\infty_B$, and this
endows $T B$ with a flat, torsion-free connection $\nabla_{GM}$ --- the
Gauss-Manin connection
of the local system $R^1f_\ast\R$. The holonomy of this connection coincides
with the monodromy of the covering $\Lambda$ (indeed, the horizontal tangent
spaces may be identified with the first homology groups of the fibres
with real coefficients).

Let $\X=R^1f_\ast\R/R^1f_\ast\Z$ be the dual family,
with projection $\hat f\colon \X\to B$. Dualizing the isomorphism
$\mbox{Vert}\,{TX}\simeq  f^\ast T^\ast B$ we get a new isomorphism
$\mbox{Vert}\,{T\X}\simeq \hat f^\ast  TB$; combining this with the splitting
of the Atiyah sequence
$$ 0 \to \mbox{Vert}\,{T\X} \to T\X \to \hat f^\ast TB \to 0 $$
provided by the Gauss-Manin connection (which can be regarded
as a connection on $T\X$), one has a splitting
$$T\X\simeq \hat f^\ast TB \oplus \hat f^\ast TB\,.$$
By letting $J(\alpha,\beta)=(-\beta,\alpha)$ this induces a complex
structure on $\X$, such
that the holomorphic tangent bundle to $\X$ is isomorphic,
as a smooth bundle, to $\hat f^\ast TB\otimes\C$.

We shall systematically use on $X$  local action-angle coordinates
$(x^1,\dots,x^g,y_1,\dots,y_g)$:   thus
the $x$'s are local coordinates on $B$, and for fixed values of the
$x$'s, the $y$'s are flat coordinates on the corresponding torus, dual to
an integral homology basis. On $\X$ we consider local
coordinates
$(x^1,\dots,x^g,w^1,\dots,w^g)$ such that the $w$'s are dual coordinates
to the $y$'s. Local holomorphic coordinates on $\X$ are given by
$z^j=x^j+iw^j$.

In this relative context it is natural to consider the fibre product
$Z=\fp$ of the fibrations $X$ and $\X$.
We shall denote by $p$, $\hat p$ the projections of   $Z$
onto its factors. On $Z$ there is a Poincar\'e bundle $\Pcc$ which may be
described in an intrinsic way, however, it is enough to say that $\Pcc$ a line
bundle on $\fp$  equipped with a  $U(1)$ connection $\nabla_\Pcc$ whose
connection form may be written in a suitable gauge as
$$\bb A=2i\,\pi\,\sum_{j=1}^gw^j\,dy_j\,.$$ Moreover, $\Pcc$
has the property that for every $\xi\in \X$, $\Pcc_{\vert\hat p^{-1}(\xi)}$
is isomorphic to $\L_\xi$ (the line bundle parametrized by $\xi$)
as a $U(1)$ bundle.

If $S$ is a closed submanifold of $X$,
we define $Z_S=S\times_B \X$, with projections $p_S$, $\hat p_S$ onto the
two factors, and denote $\Pcc_S=\Pcc_{\vert Z_S}$.
We shall assume that $\X^S=\hat p(Z_S)$ is a closed submanifold
of $\X$, and that $\hat p_S\colon Z_S\to \X^S$ is a submersion
(in the remainder of this section we shall tacitly understand that these
conditions are satisfied, while in the subsequent sections
they will hold as a consequence of other assumptions).
We  consider the exact sequence
\begin{equation}\label{GMsplits} 0 \to \hat p^{\ast}_S\Omega^1_{\X^S} \to
\Omega^1_{Z_S}
\stackrel{r}{\relbar\joinrel\to}\Omega^1_{Z_S/\X^S} \to 0\end{equation}
which defines the sheaf $ \Omega^1_{Z_S/\X^S} $ of $\hat p_S$-relative
differentials. The Gauss-Manin connection
$\nabla_{GM}$ provides a splitting of this exact sequence.
Analogously, if $\hat S$ is a closed submanifold of $\X$, we have a split
exact sequence
\begin{equation}\label{otherdiff} 0 \to p^{\ast}_{\hat S}\Omega^1_{X^{\hat
S}} \to
\Omega^1_{Z_{\hat S}}
\stackrel{\hat r}{\relbar\joinrel\to}\Omega^1_{Z_{\hat S}/X^{\hat S}} \to
0\end{equation}
which defines the sheaf $ \Omega^1_{Z_{\hat S}/X^{\hat S}} $ of $p_{\hat
S}$-relative
differentials. For every sheaf
$\E$ of $\Cc^\infty_S$-modules endowed with a flat
connection $\nabla$, one defines the following differential
operators:

(i) the operator $$\nabla^{\E}\colon
p^\ast_S\E\otimes\Pcc_S\otimes\Omega_{Z_S}^{\bullet}\to
p^\ast_S\E\otimes\Pcc_S\otimes\Omega_{Z_S}^{\bullet+1}\,,$$ obtained
by coupling the pullback of the connection
$\nabla$ with the connection of the Poincar\'e sheaf;

(ii) the operators $\nabla_r^\E$, $\nabla_{\hat r}^\E$ obtained
by composing $\nabla^\E$ with the projections $r$, $\hat r$ onto the relative
differentials.
One has $(\nabla_r^\E)^2=(\nabla_{\hat r}^\E)^2=0$.

We shall consider the higher direct images
$R^i\hat p_{S,\ast}\ker\nabla_r^\E$, which
are the cohomology sheaves of the complex
\begin{equation}\label{complex} \hat p_{S,\ast}(p^\ast_S\E\otimes\Pcc_S)
\stackrel{\nabla_r^\E}{\relbar\joinrel\to}
\hat p_{S,\ast}(p^\ast_S\E\otimes\Pcc_S\otimes \Omega^1_{Z_S/\X^S} )
 \stackrel{\nabla_r^\E}{\relbar\joinrel\to}
\hat p_{S,\ast}(p^\ast_S\E\otimes\Pcc_S\otimes \Omega^2_{Z_S/\X^S} ) \to
\dots\end{equation}

As in the usual theory of the Fourier-Mukai transform, it is
convenient to introduce a WIT notion.\footnote{Let us recall that ``WIT''
stands for ``weak index theorem.''}

\begin{defin} The pair $(\E,\nabla)$ is said to be WIT$_k$
if $R^i\hat p_{S,\ast}\ker\nabla_r^\E=0$ for $i\ne k$.
\end{defin}

Now  we want to state  a condition
for the sheaves $R^j\hat p_{S,\ast}\ker\nabla_r^\E$ to admit a connection
induced, so to say,
by the part of the operator $\nabla^\E$ complementary to $\nabla^\E_r$ .
The splitting
of the exact sequence (\ref{GMsplits}) provided by the Gauss-Manin connection
$\nabla_{GM}$ allows one to make a splitting
$$\nabla^\E=\nabla_r^\E+\hat\nabla^\E\,.$$
The $\hat\nabla^\E$ operator  induces connections on the higher direct images
$R^j\hat p_{S,\ast}\ker\nabla_r^\E$ provided it anticommutes with the
operator $\nabla_r^\E$.
The anticommutator
$\nabla_r^\E\circ\hat\nabla^\E+\hat\nabla^\E\circ\nabla_r^\E$
may be regarded as an operator $p^\ast_S\E\otimes\Pcc_S\to
p^\ast_S\E\otimes\Pcc_S\otimes\Omega^2_{Z_S}$
and as such it coincides with the restriction to $Z_S$ of $1\otimes \bb F$,
where
$\bb F$ is the curvature of the connection $\nabla_\Pcc$ of the Poincar\'e
bundle.
As a consequence, we have:
\begin{prop} \label{connec}
Assume that the sheaf $\E$ is supported on a
closed submanifold $S\subset X$,    the sheaf $R^j\hat
p_{S,\ast}\ker\nabla_r^\E$ is
supported on a closed submanifold $\hat S\subset\X$, and   the
curvature operator $\bb F$ vanishes on $S\times_B\hat S\subset Z$. Then the
operator
$\hat\nabla^\E$ induces a  connection on the sheaf $R^j\hat
p_{S,\ast}\ker\nabla_r^\E$.
\end{prop}

Eventually, we may introduce the Fourier-Mukai transform
we shall study in the remainder of this paper.

\begin{defin} If  the pair $(\E,\nabla)$ is  WIT$_k$
and satisfies the condition in Proposition \ref{connec},
the pair $(\hat\E,\hat\nabla)$, where $\hat\E=R^k\hat
p_{S,\ast}\ker\nabla_r^\E$
and $\hat\nabla$ is the connection induced as in Proposition \ref{connec},
is called the Fourier-Mukai transform of  $(\E,\nabla)$ .
\end{defin}

We end this section with an easy lemma which  is useful when checking if
the WIT
property holds for some sheaf and connection.

\begin{lemma}\label{lemmafibres} 
Let $(\E,\nabla) $ be a local system 
supported on a closed submanifold $S$ of $X$ which intersects every fibre
$X_b$ along a closed submanifold $S_b$.
For every  $j=1,\dots, g$ there is a canonical isomorphism
\begin{equation}\label{bcal}(R^j\hat p_{S,\ast} \ker\nabla_r^\E)_{\vert \X_b}
\stackrel{\sim}{\longrightarrow} R^j\hat p_{S_b,\ast} 
\ker\nabla_{1}^{\E_b} \,,\end{equation}
where $b=\hat p(\xi)$, $\hat p_b\colon X_b\times \X_b \to \X_b$ is the  
projection onto $\X_b$ and $\E_b$ is the restriction of $\E$ to $S_b$.
\end{lemma}

\begin{proof} The restriction 
$(R^j\hat p_{S,\ast} \ker\nabla_r^\E)_{\vert \X_b}$ is
defined as $\hat\jmath_b^{-1} R^j\hat p_{S,\ast} \ker\nabla_r^{\E}
\otimes_{\hat\jmath_b^{-1}\Cc^\infty_{\X}} \Cc^\infty_{\X_b}$
(here $j_b\colon S_b\to X$ and $\hat\jmath_b\colon\X_b\to\X$
are the natural inclusions).
The result is proved by applying the topological base change
\cite{Manin} to the diagram
$$\xymatrix{
S_b\times\X_b \ar[r]^{j_b\times\hat\jmath_b} \ar[d]_{\hat p_b} 
& S\times_B\X \ar[d]^{\hat p} \\
\X_b \ar[r]^{\hat\jmath_b} & \X}$$
 \end{proof}

\subsection{Fibres and Lagrangian sections \label{fibres+sections}}
In studying the transformation of local systems supported
by Lagrangian submanifolds we start by considering
the case where the submanifold is either a fibre or
a Lagrangian section.
The first case is the simplest to deal with.
It is enough to consider the case $\mbox{rank}\,{\frak L} =1$,
since the higher rank case reduces immediately to this.
The isomorphism class of the
local system ${\frak L}^\ast$ singles out a point in $\X$,
which we denote by $[{\frak L}^\ast]$. Since $X_b\times_B\X\simeq
X_b\times\X_b$, we obtain
the usual ``tautological'' property of the Fourier-Mukai transform.

\begin{prop} The pair $(\L,\nabla)\equiv\frak L$
is WIT$_g$, and the sheaf $\hat \L=R^g\hat p_\ast\ker\nabla_r^{\L}$
is isomorphic to the skyscraper $\C([\frak L^\ast])$.
\end{prop}

Now we construct a transform for $U(1)$ local systems
supported on Lagrangian sections of $X\to B$. This will generalize
the tautological correspondence that in the absolute case holds between
skyscrapers of length one on a torus and $U(1)$ local systems on the dual
torus.
The transform will produce holomorphic line bundles
on $\X$ with compatible $U(1)$ connections which satisfy some
further conditions.

Let $S\subset X$ be the image of a Lagrangian section of $X\to B$,
and ${\frak L}\equiv(\L,\nabla)$  a $U(1)$ local system on $S$.

\begin{prop} 1.  The pair $(\L,\nabla)$
is WIT$_0$;

2.   $\hat\L=\hat p_{S,\ast}\ker\nabla_r^{\L}$ is a rank-one locally free
$\Cc_\X^\infty$-module.
\end{prop}
\begin{proof} Both claims follows from Lemma \ref{lemmafibres} and the
absolute case.
\end{proof}

Since $\bb F_{\vert S\times_B\X}=0$ the conditions of
Proposition \ref{connec} are met, so that $\hat\L$ carries
a $U(1)$ connection $\hat\nabla$.
Let us express this connection in (action-angle) coordinates.  We write the local
equations  of $S$ as
$y_j=\epsilon_j(x)$; as $S$ is Lagrangian, one has
$\displaystyle\frac{\partial \epsilon_j}{\partial x^m}=
\frac{\partial \epsilon_m}{\partial x^j}$. Moreover, the $x$'s
can be thought of as local coordinates on $S$. If the connection form of
$\nabla$ is
$A=i\sum_{j=1}^gA_j(x)\,dx^j$, with
$\frac{\partial A_j}{\partial x^\ell}=\frac{\partial A_\ell}{\partial x^j}$,
then $\hat\nabla$ may be
represented by the connection form
$$\hat A = i\sum_{j=1}^gA_j(x)\,dx^j- 2i\pi\sum_{j=1}^g
\epsilon_j(x)\,dw^j\,. $$
In these coordinates the components of the connection form $\hat A$
do not depend on the $w$'s. Moreover, both the horizontal
 and vertical part (with respect to the splitting given by the Gauss-Manin
connection)
are flat, and in particular, the restriction of $\hat\nabla$ to any fibre
$\X_b$ of $\X\to B$ is flat.

\begin{remark}
The independence of the components $\hat A$
on the $w$'s   can be stated invariantly in a variety of ways.
For instance, one can use the fact that the zero-section of $\X$ makes the
latter into
a (trivial) principal $T^g$-bundle over $B$; then,  $\hat\nabla$ commutes
with the
action of $T^g$ on $\X$.
\label{winvar}\end{remark}
The Hodge components of curvature
form $\hat F$ of this connection may be written ---
recalling that in the complex structure we have given to $\X$ the
coordinates  $z^j=x^j+iw^j$ are complex holomorphic  --- as
$$\hat F^{2,0}=\frac{\pi}{2} \sum_{k,j} \frac{\partial \epsilon_j}{\partial
  x^k}\, dz^k\wedge dz^j$$
$$\hat F^{0,2}=-\frac{\pi}{2} \sum_{k,j} \frac{\partial \epsilon_j}{\partial
  x^k}\, d\bar{z}^k\wedge d\bar{z}^j$$
$$\hat F^{1,1}=\frac{\pi}{2} \sum_{k,j} (\frac{\partial \epsilon_k}{\partial
  x^j}+\frac{\partial \epsilon_j}{\partial
  x^k}) \,dz^k\wedge d\bar{z}^j.$$
Since $S$ is Lagrangian we have  $\hat F^{0,2}=\hat F^{2,0}=0$, so that
$\hat\L$ may be given a holomorphic structure
compatible with the connection $\hat\nabla$. Moreover,
we have
$$
\hat F^{1,1}=\pi \sum_{k,j} \frac{\partial \epsilon_k}{\partial
  x^j} dz^k\wedge d\bar{z}^j.
$$

\begin{defin} The Fourier transform of
$(S, {\frak L} )$ is the pair
$(\hat \L,\hat \nabla)$.
\end{defin}

\subsection{The non-transversal case \label{skew}}
The results in section \ref{fibres+sections} can be generalized to local
systems  supported on Lagrangian submanifolds of $X$ other than sections. This
allows us to enlarge  the ``dual'' category on which the inverse Fourier-Mukai
transform  (defined in section \ref{inverse}) acts, to a category of
sheaves with connection  (satisfying some suitable conditions) supported by
complex submanifolds of $\X$. Such sheaves arise naturally in Fukaya's
treatment of mirror symmetry   \cite{Fu}. We shall be interested in
transforming local systems  supported on a submanifold $S$ of $X$ such that:

\begin{itemize}
\item[(C1)] $S$ is a Lagrangian subvariety of $X$;
\item[(C2)] the intersection $S_b = S\cap X_b$ of $S$ with a fibre of $X$,
when nonempty,
is a   (possibly affine) subtorus $ S_b$ of $X_b$ whose dimension does not
depend on $b$.
\end{itemize}
Let $\frak L\equiv(\L,\nabla)$  be a $U(1)$ local system on $S$.
Its Fourier-Mukai transform
at the sheaf level is
$$\hat \L =R^{m}\hat p_{S,\ast} \ker \nabla_r^{\L} $$
where $m$ is the dimension of the tori $S_b$. Indeed, one has:

\begin{prop} If the  conditions C1
and C2 are satisfied, the sheaf $\cal L$ is WIT$_{m}$.
\end{prop}

\begin{proof} It follows from Lemma \ref{lemmafibres} and Proposition
\ref{absoluteprop}. \end{proof}

 Lemma \ref{lemmafibres} and Proposition \ref{absoluteprop} also imply that
after restriction to its support, $\hat \L$ is a line bundle. We shall now
show that, under some suitable conditions on
the support $S$,  the transform   $\hat\L$ is   supported on a complex
submanifold $\hat S$
of the dual family $\X$.
More precisely, we  assume:

\begin{itemize}\item[(C3)] the vertical tangent spaces of the family of
subtori $\{S_b\}_{b\in f(S)}$ are parallelly transported by the Gauss-Manin
connection
$\nabla_{GM}$    regarded as a connection in $TX$.
\end{itemize}

This   requirement can  be translated into a more
explicit form in terms of the action-angle coordinates $(x,y)$ we have
previously introduced, in that it amounts to the condition that the family
of subtori $\{S_b\}$ can be
  written as
\begin{displaymath}
\sum_{j=1}^{g} a_i^j \,y_j+\chi_i=0\,, \qquad
 i=1,\dots ,g-m
\end{displaymath}
with the matrix $a^j_i$ constant and the $\chi_i$'s  local functions on $B$.
\begin{lemma} Conditions C1, C2 and C3 imply that $f(S)$ is a submanifold
of $B$ of dimension
$k=g-m$, and that it can be parametrized by the first $k$  action
coordinates $x^j$.
\end{lemma}
\begin{proof} The first claim follows from the fact that the horizontal part
of the tangent space to $S$ has constant dimension; the second from
the Lagrangian condition which implies that the local equations of $f(S)$
in $B$
are linear in the action coordinates.
\end{proof}

\begin{prop}\label{transf}
Let $(S, \L,\nabla )$  be a local system supported on a Lagrangian submanifold
$S$   fulfilling the conditions C1 and C2. The  condition C3 is satisfied
if and only if the support $\hat S$ of
the transform $\hat \L$ is a complex submanifold of $\X$.
\end{prop}

\begin{proof} A proof is given in Appendix B. \end{proof}

\begin{remark} In our setting there is no constraint on the dimension
of $X$, the latter space is assumed to be just symplectic, and
we consider local systems supported on Lagrangian submanifolds of $X$.
On the other hand, string-theoretic mirror symmetry assumes,
on physical grounds, that $X$ is a (usually 3-dimensional) Calabi-Yau
manifold, and one considers special Lagrangian supports.\footnote{Let us
recall that a special Lagrangian submanifold of a Calabi-Yau $n$-fold
$X$ is an oriented real $n$-dimensional submanifold $Y$ which is Lagrangian
w.r.t.~the K\"ahler form of $X$, and such that one can choose a global
trivialization of the canonical bundle of $X$ whose imaginary
part vanishes on $Y$. For more details cf.~\cite{HL}.}
In this case, the condition that $S$ is special Lagrangian implies, for
$k=1$,  that the
coefficients
$a_i^j$ are constant, so that this is a particular case within our treatment.
On the contrary, for $k=2$ the speciality  property seems to be unrelated
to the conditions that ensure the support $\hat S$ to be complex holomorphic.
\end{remark}

\begin{prop} Under the conditions of Proposition \ref{transf},
the operator $\hat\nabla^\L$ (cf.~section \ref{framework}) induces on
$\hat\L$   a $U(1)$ connection.
\end{prop}
\begin{proof} This will use the proof of Proposition \ref{transf} given in
Appendix B.
We know that  $\hat\nabla^\L$
induces a connection on the Fourier-Mukai transform if the curvature $\bb
F$ of the
Poincar\'e bundle on $Z=\fp$  vanishes on $S\times_B\hat S$, where
$S$ and $\hat S$ are the supports of $\L$ and $\hat \L$, respectively.
In view of the form of $\bb F$, this condition is met if for each $b\in B$ the
intersections of $S$ and $\hat S$ with the fibres $X_b$, $\X_b$ yield
subtori of $X_b$, $\X_b$ that are normal to each other. But looking at the
equations
of the supports, (\ref{sist1}) and (\ref{supp2}), and comparing with the
absolute case (Proposition \ref{absoluteprop}), we see that this condition is
fulfilled. \end{proof}

 We shall now prove that $\hat \L$, as a line bundle on $\hat S$, has
a holomorphic structure.
Let $\hat \nabla$ be the connection induced on $\hat \L$.

\begin{prop}\label{holotrans}
If the support $\hat S$ of the transformed sheaf $\hat \L$ is
a complex submanifold of $\X$, then   $\hat\nabla$ induces a holomorphic
structure on $\hat \L$.
\end{prop}
\begin{proof}
The connection 1-form of the connection $\nabla $ can be written in an
appropriate gauge as
$$A = i   \sum_{j=1}^k \alpha_j (x^1, \dots ,x^k)\,dx^j + 2i\,\pi
\sum_{\ell=1}^{g-k} \xi^\ell\,
dy_\ell\,,$$   with the quantities $\xi^\ell$
constant. From the proof of Proposition \ref{absoluteprop} given in the
Appendix A  we know that the transformed connection $\hat \nabla$ is given in
local  action-angle coordinates by the 1-form\footnote{Also in this case we find that the
transformed connection $\hat\nabla$ satisfies the condition of Remark
\ref{winvar}.} $$\hat{A}=-2i\,\pi\sum_{\ell=g-k+1}^g \chi_{\ell}(x^1,\dots
,x^k)\,dw^{\ell}   +  i\sum_{j=1}^k \alpha_j (x^1, \dots ,x^k)\,dx^j\,.$$
Rewriting this in terms of $w^1,\dots ,w^k$ we obtain
$$
\hat A=-2i\,\pi\sum_{\ell=g-k+1}^g\sum_{j=1}^k\chi_\ell(x^1,\dots,x^k)\,\tilde\gamma^\ell_j\,dw^j
+i\sum_{j=1}^k \alpha_j (x^1, \dots ,x^k)\,dx^j$$
where
the  coefficients $\tilde\gamma^\ell_j$ are constant.
Since $d(\sum_j\alpha_j\,dx^j)=0$   because of the flatness of $\nabla$, it
follows that
the curvature of $\hat \nabla $ is given by
$$
\hat F=-2i\,\pi\,\sum_{\ell=g-k+1}^g \sum_{j,m=1}^k \frac{\partial
\chi_{\ell}}{\partial
x^j}\,\tilde\gamma^\ell_m\,dx^j\wedge dw^m\,.
$$
Since  $\displaystyle\tilde\gamma^{g-k+j}_m=\frac{\partial
\zeta^{g-k+j}}{\partial x^m}$,
where the functions $\zeta^{g+j}$ are those of the equations (\ref{sist1}),
the condition $\hat F^{0,2}=0$ can be written as
$$
\sum_{\ell=g-k+1}^g\left[\frac{\partial \zeta^{\ell}}{\partial x^j}
\frac{\partial \chi_\ell}{\partial x^m}-
\frac{\partial \zeta^{\ell}}{\partial x^m}
\frac{\partial \chi_\ell}{\partial x^j}\right]=0\,, \qquad 1\leq j<m\leq k.
$$
But this is the system of equations (\ref{eq4}), therefore when $S$
is Lagrangian, this condition is automatically satisfied.
\end{proof}

\begin{remark} (The higher rank case.) So far  we have
for simplicity considered only the transformation of local systems of rank one.
However the higher rank case, under the same conditions, can be treated along
the same lines, obtaining on the $\X$ side holomorphic vector bundles of the
corresponding rank supported on complex submanifolds of $\X$. \end{remark}

\subsection{Invertibility\label{inverse}}
In this section we shall prove that the Fourier-Mukai transform we have
defined inverts.
However we shall only  discuss the inverse transform of rank 1 sheaves.
The higher rank case requires to consider Lagrangian submanifolds
of $X$ which ramify over $B$, and this will be done in a future paper.

We shall therefore consider a holomorphic line bundle $\hat \L$
supported on a $k$-dimensional complex submanifold $\hat S$ of $\X$, equipped
with a compatible $U(1)$ connection $\hat\nabla$. Moreover, we shall assume
that:

\begin{itemize}
\item[(D1)] $\hat S$ intersects
the fibres of $\X $ along affine subtori of complex dimension $k$;
\item[(D2)] the horizontal part of the connection $\hat\nabla$ is flat
(horizontality  is
given by the Gauss-Manin connection);
\item[(D3)] the connection $\hat\nabla$ is invariant under the action of
$T^g$ on $\X$
(cf.~Remark \ref{winvar}).
\end{itemize}
These conditions allow us to write the local connection form of $\hat\nabla$ as
$$ \hat A=i\,\sum_{j=1}^k
\alpha_j(x^1,\dots,x^k)\,dx^j+2i\,\pi\sum_{j=1}^k\beta_j(x^1,
\dots,x^k)\,dw^j \,, $$ where the functions $\alpha_j$ satisfy (as a
consequence of D2) the
closure condition
$\frac{\partial\alpha_j}{\partial
x^\ell}=\frac{\partial\alpha_\ell}{\partial x^j}$.
This shows that the restriction of $\hat\nabla$ to any fiber $\X_b$ of
$\X\to B$
yields a flat connection on $\hat\L_{\vert\X_b}$.

Let $p_{\hat S},\hat p_{\hat S}$ the canonical projections of $X\times_B
\hat S$ onto its
factors. We consider the operator
$$\nabla_{\hat r}^{\hat\L}=\hat r\circ(\hat p^\ast_{\hat S}\nabla^{\cal L}
\otimes
1+1\otimes \nabla_{\Pcc^\vee_{\hat S}})$$ and in terms of it we define a
Fourier-Mukai
transform from sheaves on $\X$ to sheaves on $X$ (notice that we twist with
the dual
Poincar\'e bundle $\Pcc^\vee$).
\begin{prop} $\cal L$ is WIT$_k$, and
$\L=R^k p_{\hat S,\ast}\ker\nabla_{\hat r}^{\hat\L}$ is  supported on a
Lagrangian submanifold $S$ of $X$ such that every
intersection $S_b = S\cap X_b$
is an affine subtorus of $X_b$ of dimension $g-k$ (when nonempty). Moreover
the family
of subtori $S_b$ is parallelly
transported by the Gauss-Manin connection $\nabla_{GM}$. Finally,
a flat connection $\nabla$ is naturally induced on $\L$.
\end{prop}

\begin{proof} The WIT condition follows immediately from Lemma
\ref{lemmafibres}. To show  the remaining part of the claim  we write  local
equations  for $\hat S$  as
\begin{displaymath}
 \left\{\begin{array}{ll}
x^{k+j} = \zeta^{k+j}(x^1, \dots ,x^k)\,,&   j=1,\dots , g-k\\[3pt]
w^{k+j} = \sum_{i=1}^{k} P_i^{k+j} (x^1, \dots , x^k)\,w^i + Q^{k+j}(x^1,
\dots , x^k)\,,&  j=1,\dots ,
g-k.
 \end{array}\right.
\end{displaymath}
Performing a fibrewise transform we obtain the following equations
for the support $S$ of the transform $\L$:
\begin{displaymath}
 y_l + \sum_{m=k+1}^gP^{m}_l (x^1,\dots ,x^k)\, y_{m}
+\beta_l (x^1,\dots ,x^k) = 0
\end{displaymath}
where $l=1, \dots ,k$.  It remains to show that $S$ is Lagrangian and that the
family $\{S_b\}_{b\in \hat f(S)}$ is parallelly transported by the
Gauss-Manin connection. The latter point follows from  the complex
structure of
$\hat S$ (cf.~Proposition \ref{transf}): the Cauchy-Riemann equations for
$\hat S$  imply
that  the coefficients $P_l^{k+j}$ and $Q^{k+j}$ are constant. As far as the
Lagrangian property of $S$ is concerned, the holomorphicity of $\hat S$ and
$\hat\L$ imply
the  equations (\ref{sist2}) in the proof of Proposition \ref{transf} (in
Appendix B).
Therefore $S $ is Lagrangian.  Observe that the transformed connection $\
\nabla $ has a
1-form given by
$$  A=i\sum_{j=1}^k \alpha_j (x^1,\dots ,x^k)\, dx^j  -
2i\,\pi\sum_{m=k+1}^gQ^{m}\,dy_{m} \,,$$
whence we can immediatly deduce its flatness.
\end{proof}

\section{Conclusions\label{Conclusions}}
We summarize here the main results of this paper.
We have shown that a suitably defined Fourier-Mukai transform $\F$ maps a
$U(1)$ local
system supported on a Lagrangian subvariety $S$ of $X$ satisfying  the
conditions
C1, C2 and C3 (cf.~section \ref{skew}) into a holomorphic line bundle $\hat
\L$  supported
a complex  subvariety $\hat S$ of $\hat X$; moreover   $\hat \L$ is
endowed with a $U(1)$ connection such that conditions D1, D2, D3 (section
\ref{inverse}) are satisfied.

Conversely, if we start with a  holomorphic line bundle supported on a
complex subvariety $\hat S$ of $\hat X$ equipped with a $U(1)$ connection
$\hat \nabla $ such that conditions D1, D2, D3 are satisfied,  we define
a dual  Fourier-Mukai transform $\hat \F$ that maps such objects into a
$U(1)$ local system supported on a Lagrangian subvariety $S$ such that
conditions C1, C2 and C3 are fulfilled. The explicit forms of the two
transforms
we have written in sections \ref{fibres+sections} and \ref{skew}
show that the transforms are one the inverse of the other.
This  parallels the classical  result in \cite{M2} and generalizes 
the one in \cite{AP}, whose authors consider the
case where $X$ and $\hat X$ are
$S^1$-fibrations  over $S^1$ ($\hat X$ is actually an elliptic curve) and
$\L$ is a local
system on an affine line $S\subset X$. Observe that in this case the
conditions C1, C2, C3  and  D1, D2, D3 are trivially satisfied.

Finally, we would like to comment upon the relation of the construction we
have described
in this paper with Fukaya's homological mirror symmetry.
First we notice that, in the absence of the B-field and with no singular
fibres,
our ``mirror manifold'' $\X$ coincides with Fukaya's, also taking into account
its complex structure. Let $S$ be a Lagrangian submanifold of $X$, and
$\beta=(\L,\nabla)$
a local system on it. Fukaya proposes to construct on $\X$ a coherent sheaf
whose fibre at a point $(b,\alpha)\in\X$ (where
$\alpha=(L_\alpha,\nabla_\alpha)$
is a local system on the fibre
$X_b$) is given by the Floer homology
$$HF^\bullet((X_b,\alpha),(S,\beta))\,.$$
This homology may be proved \cite{FuFlo} to be isomorphic to
$$H^{\bullet-\eta(X_b,S)}(S\cap X_b, \cal Hom_{\nabla}(\L_\alpha,\L))\,,$$
where $\eta(X_b,S)$ is a Maslov index, and  $\cal Hom_{\nabla}(\L_\alpha,\L)$
is the sheaf of $\nabla$-compatible morphisms between $\L_\alpha$ and $\L$.
It is not difficult to show that only one of these cohomology
groups does not vanish (in the correct degree), and that it is
isomorphic,  up to a dual,  to
the fibre of our transform $\hat\L$. However, the concrete construction done in
\cite{Fu} is not in terms of Floer homology, but it is an \emph{ad hoc} one,
which may be compared with ours when $X=T^{2g}$, $B=T^g$ and $S$ is a
Lagrangian embedding of $T^g$. In this case the vector bundle constructed on
$\X$ coincides with ours.

It should be noted that our construction provides on the ``mirror side'' $\X$
more data, in that we obtain on $\hat\L$ a connection.
 It is interesting to note that this connection
is not invariant under Hamiltonian diffeomorphisms of $X$, while the
remaining geometric data on $\X$ are.

\section*{Appendix A}
We provide here a sketch of the proof of Proposition \ref{absoluteprop}.
It involves a number of computations but it is conceptually very easy, being
a generalization of the proof of Proposition \ref{prop1}.
For further details we refer to \cite{th}.

We first consider the case when $S$ is a 1-dimensional
affine subtorus of $T$. The direct image $R^i\hat{p}_{S,\ast} \ker
\nabla_1^{\cal L}$ is by definition the sheaf associated to
the presheaf
$\hat{U}\rightsquigarrow
H^i(S\times\hat{U},\ker\nabla_1^{\cal L})\simeq
H^i\left(\Omega^{\bullet,0}(p^\ast_S
L\otimes\Pcc_S)(S\times\hat{U}),\nabla_1^{\cal L}\right)$.
When $i=0$, take an element $s$ in
$H^0(S\times\hat{U},\ker\nabla_1^{\cal L})$ and consider its
restriction to $S\times\{y\}$, with $y\in\hat{U}$, which is a
global section of ${\cal L}\otimes \Pcc_{\vert S\times\{y\}}$.
This means that $\Pcc_{\vert S\times\{y\}}$ is non trivial for
every $y$ in the complement of $\ker\psi$, which is a dense
subset of $\T$ (here $\psi$ is the natural map
$\psi:\T\rightarrow
Hom(\pi_1(S),U(1))$).
 The same holds for ${\cal L}\otimes
\Pcc_{\vert S\times\{y\}}$. Since $\nabla_1^{{\cal L}}
s_{\vert S\times\{y\}}=0$
for every $y\in\T$, the
restriction of $s$ to $S\times\{y\}$ vanishes for $y$ in a
dense subset of $\T$, so that $s$ vanishes everywhere.

When $i=1$ we need to write  the equation of $S$ explicitly.
For simplicity we only give some details in the case   $\dim T=2$. Let
$(y^1,y^2)$ be flat coordinates on $T$ and   $(w_1,w_2)$ flat
dual coordinates on $\T$. We pick a gauge where the Poincar\'e bundle has
an automorphy
factor
$$a_\Pcc(y^1,y^2,w_1,w_2,\lambda^1,\lambda^2,\mu_1,\mu_2)=
e^{2\pi i(\lambda^1\,w_1+\lambda^2\,w_2)}.$$
The equation of $S$ in the universal cover of $T$ is given
by an affine line $y^2=ay^1+b$.
Let $A=\bar{y}_1\,dy^1$
be the connection form of the local
system $({\cal L},\nabla^{{\cal L}})$ on $S$. We need to
compute $H^1(S\times\hat{U},\ker\nabla_1^{\cal L})$. So
take an element
$\tau\in(\Omega^{1,0}(p^\ast_S {\cal
L}\otimes\Pcc_S)(S\times\hat{U}),\nabla_1^{{\cal L}})$.
Observe that $\tau$ is closed with respect to $\nabla_1^{\cal
L}$ because $\dim S=1$. If we let
$$\tau=\phi(\xi,w_1,w_2)\,d\xi\,$$
where $\xi$ is the natural coordinate on $S$,
the automorphy condition satisfied
by $\tau$ can be expressed  in the form
$$\phi(\xi+\sqrt{p^2+q^2},w_1,w_2)=e^{p(w_1+\bar{w_1})+qw_2}\phi(\xi,w_1,w_2)$$
having set $a=q/p$ with $q,p$ coprime.

Suppose that $\tau$ is exact so that we can write
$\tau=\nabla_1^{\cal L} s$ where $s\in
C^\infty(S\times\hat{U},\ker\nabla_1^{\cal L})$. Then $s$ can be
written in the form
$$s(\xi,w_1,w_2)=\int_0^{\xi}\phi(u,w_1,w_2)du+c(w_1,w_2)$$
but this is well defined if and only if the automorphy condition
 is satisfied, and one can easily check that this amounts to
$$
c(w_1,w_2)(1-e^{2\pi
i(p(w_1+\bar{w_1})+qw_2)})=
-\int_0^{\sqrt{p^2+q^2}}\,\phi(u,w_1,w_2)\,du.
$$
This equation may be solved for $c$
in the complement of the set $\hat{S}$
defined by
$$w_2=-\frac{1}{a}w_1-\frac{\bar{w_1}}{a};$$
thus, arguing as in Proposition \ref{prop1}, we obtain that the
support of $R^1\hat{p}_{S,\ast} \ker \nabla_1^{\cal L}$ is exactly
$\hat{S}$.

To compute the sheaf $R^1\hat{p}_{S,\ast} \ker \nabla_1^{\cal L}$ we note
that the map
\begin{eqnarray*} \varpi\colon \Omega^{1,0}(p^\ast_S
{\cal L}\otimes\Pcc)(S\times\hat{U}) &\to& \Cc^\infty(\hat{S}\cap\hat{U}) \\
\tau &\mapsto&  -\int_0^{\sqrt{p^2+q^2}}\,\phi(u,w_1,w_2)\,du
\end{eqnarray*}
is surjective: if
$f\in\Cc^\infty(\hat{S}\cap\hat{U})$ and $s$ is a
section of the Poincar\'e bundle over
$S\times\hat{S}$,
then the 1-form
$\tau=\phi\,d\xi$ defined by
$$\phi(\xi,w_1,w_2)=\beta \,s(\xi,w_1)\,f(w_2),$$
with $1/\beta=-\int^{\sqrt{p^2+q^2}}_0s(u,0)\,du$,
satisfies $\varpi(\tau)=f$ and the correct automorphy condition.
Thus,  $H^1(S\times\hat{U},\ker\nabla_1^{\cal L})=
\Cc^\infty(\hat{S}\cap\hat{U})$.

The transformed sheaf is endowed with a flat connection
induced by $\nabla_2^{{\cal L}}$, the $(0,1)$ part of
the connection $p^*_S\nabla^{{\cal L}}\otimes1+1\otimes\nabla_{\Pcc_S}$,
because the anticommutator between $\nabla_1^\L$ and $\nabla_2^\L$ equals
the curvature of the Poincar\'e bundle, which  restricts to zero on
$S\times\hat{S}$
(cf.~Proposition \ref{connec}).

Of course, $R^i\hat p_{S,\ast} \ker\nabla_1^{\cal L}=0$
for $i>1$ because $S$ is 1-dimensional.
This proof is extended to the case  $\dim S  > 1$ by using a K\"unneth
formula.

\section*{Appendix B}
Here we prove Proposition \ref{transf}. For notational convenience we
suppose that $k\leq g/2$; the complementary case $ k>g/2$
can be treated  similarly. In the  action-angle coordinates $x$, $y$
we can write the local equations for $S$ as
\begin{equation}\label{sist1}
\left\{ \begin{array}{ll}
y_{g-k+j} =\eta_{g-k+j}(x^1,\dots ,x^k,y_1,\dots ,y_{g-k})\,, &
j=1,\dots ,k\\
x^{k+i} =\zeta^{k+i}(x^1,\dots ,x^k)\,, & i=1,\dots ,g-k
\end{array}\right.
\end{equation}
Since $S$ is Lagrangian one has
\begin{equation}\label{sist2}
\left\{ \begin{array}{ll} \displaystyle
\delta_j^m+\sum_{\ell =g-k+1}^g \frac{\partial \zeta^\ell }{\partial x^j}
\frac{\partial \eta_{\ell }}{\partial y_m}=0\,,& j,m=1, \dots ,k\\[7pt]
\displaystyle
\frac{\partial \zeta^{k+i}}{\partial x^m}+
\sum_{\ell =g-k+1}^g
\frac{\partial \zeta^{\ell }}{\partial x^m}
\frac{\partial \eta_{\ell }}{\partial y_{k+i}}=0\,, & i= 1,\dots ,g-2k;m=1,
\dots ,k;
\\[7pt] \displaystyle
\sum_{\ell =g-k+1}^g\left[\frac{\partial \zeta^{\ell }}{\partial x^j}
\frac{\partial \eta_{\ell }}{\partial x^m}-
\frac{\partial \zeta^{\ell }}{\partial x^m}
\frac{\partial \eta_{\ell }}{\partial x^j}\right] =0\,, & 1\leq j<m\leq k.
\end{array}\right.
\end{equation}
The equations of the subtori $S_b$ can be written in
a linear form
\begin{equation}\label{torilineari}
y_{g-k+j}=\sum_{m=1}^{g-k}a_{g-k+j}^m(x^1,\dots ,x^k)\,
y_m+\chi_{g-k+j}(x^1,\dots ,x^k)\,, \quad j=1,\dots ,k\,.
\end{equation}
To find the equations of $\hat S$ we shall perform a fibrewise transform
and use the K\"unneth formula as in Proposition \ref{prop2}. First we split
every subtorus $S_b$ as a product of 1-dimensional tori $r_i (b)$ which
have linear equations given by
\begin{displaymath} \left\{ \begin{array}{ll}
y_l =0\,, &  l=1,\dots ,g-k,\ \ell\ne i;  \\
y_{g-k+j} = a_{g-k+j}^i (x^1,\dots ,x^k) y_i+\chi_{g-k+j} (x^1,\dots ,x^k)\,,&
  j=1,\dots ,k.
\end{array}\right.
\end{displaymath}
Observe that we can also split the local system $\L$ on $S_b$
as a box product
 of local systems $\L_i (b)$ on $r_i (b)$ where
$i=1,\dots ,g-k$. Transforming the local system $\L_i (b)$ on $r_i (b)$ we
get the
following  equation for the support of $\L_i (b)$ (see Appendix A):
$$w^i+  \sum_{\ell=g-k+1}^g\gamma_{\ell}^i(x^1,\dots ,x^k)\, w^{\ell}+ \xi^i$$
where the constant term $\xi^i$ describes the automorphy of $\L_i$ (here $i$
is fixed), and the matrix $\gamma^i_\ell$ satisfies the condition
$\sum_{j=1}^g\gamma^j_\ell\,a^i_j=0$.
Then  $\hat
S$ is the intersection of the supports $\hat r_i $, so that its equations are
of the form \begin{equation} \label{supp2}
w^{k+i}=\sum_{j=1}^k\tilde\gamma^{k+i}_j (x^1,\dots ,x^k)\,
w^j+\varsigma^{k+i}(x^1,\dots ,x^k) \,,\quad i=1,\dots , g-k\end{equation}
together with the second set of equations (\ref{sist1}). Here we have solved
with respect to $w^1,\dots ,w^k$. These equations may be used to replace the
functions $\eta$ in (\ref{sist1}), thus getting
\begin{equation}\label{sist3}
\left\{\begin{array}{ll} \displaystyle
\delta_j^m+\sum_{\ell=g-k+1}^g\frac{\partial \zeta^{\ell}}{\partial x^j}\,
a_{\ell}^m =0\,, &  j,m=1,\dots ,k\\[5pt] \displaystyle
\frac{\partial \zeta^{k+i}}{\partial x^m}+ \sum_{\ell=g-k+1}^g
\frac{\partial \zeta^{\ell}}{\partial x^m}
a_{\ell}^{k+i} =0 \,,&  i=1,\dots ,g-2k; m=1,\dots ,k
\end{array}\right.
\end{equation}
\begin{equation}\label{eq4}\sum_{\ell=g-k+1}^g\left[
\frac{\partial \zeta^{\ell}}{\partial x^j}
\frac{\partial \chi_{\ell}}{\partial x^m}-
\frac{\partial \zeta^{\ell}}{\partial x^m}
\frac{\partial \chi_{\ell}}{\partial x^i}\right]
=0\,, \quad 1\leq j<m\leq k.
\end{equation}
The  solution of (\ref{sist3}) is
\begin{equation}\label{gamma}
\frac{\partial \zeta^{k+i}}{\partial x^j}= \tilde\gamma^{k+i}_j\,, \qquad
\qquad  j=1,\dots ,k, \quad    i=1,\dots ,g-k\,.
\end{equation}
If the submanifold $S$ is Lagrangian, the conditions (\ref{gamma}) admit
solutions in $\zeta$. We must check that the support $\hat S$ is holomorphic,
i.e.,  the equations that define it fulfil the Cauchy-Riemann conditions.
The latters are satisfied if and only if the coefficients
$\tilde\gamma^{k+i}_j$
do not depend on the $x$'s, but this is true if and only if the coefficients
$\gamma_{g-k+1}^{k+i}$ are in turn independent of the $x$'s. As a result,
we have
proved that when $S$ is Lagrangian, the  tangent spaces to the $S_b$'s
are parallelly transported by $\nabla_{GM}$ if and only if $\hat S$ is
holomorphic.

One may note that the coefficients $\chi^j$ play no role in the
specification of the
complex structure of $\hat S$. Moreover, let us remark that equation
(\ref{supp2}) shows
that the intersections of the support $\hat S$ with the fibres $\X_b$ are
affine subtori.

\end{document}